\documentclass[a4paper,12pt]{article}
\usepackage{amssymb}
\usepackage{amsmath}
\usepackage{graphicx}
\usepackage{epstopdf}

\newcommand{\re}[1]{(\ref{#1})}
\def\mp{{{\cal M}^+_{\lambda,\Lambda}}}
\def\mm{{{\cal M}^-_{\lambda,\Lambda}}}

\def\Lpl{\mathcal{L}^+}
\def\Lmin{\mathcal{L}^-}
\def\Lpm{\mathcal{L}^\pm}

\newcommand{\eneqa}{\end{eqnarray}}
\newcommand{\begeqaet}{\begin{eqnarray*}}
\newcommand{\eneqaet}{\end{eqnarray*}}
\newcommand{\be}{\begin{equation}}
\newcommand{\ee}{\end{equation}}

\newcommand{\rn}{\rbig^n}

\newcommand{\rbig}{{\mathbb{R}}}

\newcommand{\lpr}{{L^\prime}}
\newcommand{\ed}{\end{document}}

\newcommand{\ep}{\varepsilon}

\def\beq{\begin{equation}}
\def\eeq{\end{equation}}

\def\Box{\hfill\framebox(0.25,0.25){}}

\newtheorem{thm}{Theorem}[section]
\newenvironment{thmbis}[1]
  {\renewcommand{\thethm}{\ref{#1}$'$}%
   \addtocounter{thm}{-1}%
   \begin{thm}}
  {\end{thm}}

\newtheorem{lem}{Lemma}[section]
\newtheorem{cor}{Corollary}[section]
\newcommand{\begeqa}{\begin{eqnarray}}

\begin{document}

\begin{center}{\bf\Large Boundary Harnack estimates and quantitative strong maximum principles for uniformly elliptic PDE}
\end{center}\smallskip

 \begin{center}
Boyan SIRAKOV\footnote{e-mail : bsirakov@mat.puc-rio.fr}\\
PUC-Rio, Departamento de Matematica,\\ Gavea, Rio de Janeiro - CEP 22451-900, BRAZIL
\end{center}\bigskip

{\small \noindent{\bf Abstract}. We give full boundary extensions to two fundamental estimates in the theory of elliptic PDE,  the weak Harnack inequality and the quantitative strong maximum principle, for uniformly elliptic equations in non-divergence form. }
%

\section{Introduction and Main Results}

In this paper we establish boundary and global versions of some important estimates in the theory of uniformly elliptic PDE in non-divergence form, the (weak) Harnack inequality and the quantitative strong maximum principle.

We consider viscosity solutions of inequalities in the form
\begin{equation}\label{first}
\Lmin[u]:=\mm(D^2u)  - b(x) |Du|\le f(x)
\end{equation}
(or $\Lpl[u]:= \mp(D^2u) + b(x) |Du|\ge -f(x)$), where $0<\lambda\le \Lambda$, $b,f$ belong to suitable Lebesgue spaces, and $\mp$ is the extremal Pucci operator
$$\mathcal{M}^-(M)=\Lambda \underset{\mu_i<0}{\sum} \mu_i+\lambda \underset{\mu_i>0}{\sum} \mu_i = \inf_{\lambda I \le A \le \Lambda I} \mathrm{tr}(AM),\quad\mathcal{M}^+(M)= -\mathcal{M}^-(-M),$$
for any symmetric matrix  $M\in \mathcal{S}_n$, where  $(\mu_i)_{i=1..n}$ are the eigenvalues of $M$. The operator $\Lmin$ can be viewed as the infimum of linear operators with fixed bounds for the coefficients.

A basic property of linear and some nonlinear uniformly elliptic PDE is the {\it strong maximum principle} (SMP). It states that if $u$ is a solution of $\Lmin[u]\le 0$ in some domain $\Omega\subset\rn$, such that $u\ge 0$ in $\Omega$ and $u(x_0)=0$ for some $x_0\in\Omega$, then $u\equiv0$ in $\Omega$. The SMP can be obtained as a consequence of {\it Hopf's lemma}, which says that if $u$ is a solution of $\Lmin[u]\le 0$ in the unit ball $B_1\subset\rn$, such that $u> 0$ in $B_1$ and $u(x_0)=0$ for some $x_0\in\partial B_1$, then $\liminf_{t\searrow0} t^{-1}u((1-t)x_0)>0$.

In the following $\Omega$ denotes a bounded $C^{1,1}$-domain in $\rn$. We will proceed under the minimal integrability requirements on the coefficients $b$ and $f$ which ensure the validity of the Hopf lemma for $\Lmin[u]\le 0$, and the solvability of $\Lmin[u]=f$ in $\Omega$, $u=0$ on $\partial \Omega$,  together with the finiteness of the quantity $u(x)/\mathrm{dist}(x,\partial\Omega)$ for any solution of this problem (we stress however that to our knowledge Theorems \ref{bqsmp}-\ref{bhi} below have not appeared before even for linear equations with bounded or smooth coefficients). Specifically, in the following we assume that for some $q>n$ and $q\ge p>p_0$ we have
 $$
 b\in L^q(\Omega), \;\; b\ge0\mbox{ in }\Omega,\qquad f\in \lpr(\Omega):= L^p(\Omega)\cap L^q(\Omega_{d_0}),
 $$
 where $\Omega_{d_0}:=\{x\in\Omega\,:\, \mathrm{dist}(x,\partial\Omega)<d_0\}$, for some fixed (small) $d_0>0$, and  $n/2<p_0=p_0(n, \lambda, \Lambda)<n$ is the optimal exponent for the validity of the ABP estimate for the Pucci extremal operators (see \cite[Theorem 9.1]{GT}, as well as \cite{FS}, \cite{E}, \cite{GS}).\smallskip

A far-reaching and well-known extension of the SMP is the {\it interior weak Harnack inequality} (WHI), a pivotal result in the regularity theory for elliptic PDE, which goes back to de Giorgi (for divergence form equations) and Krylov-Safonov (for equations in non-divergence form); see \cite[Chapters 8, 9]{GT}, \cite{Kbook}, \cite{CC}, and \cite{KS2} for the most general results. The WHI states that if $\Lmin u\le0$ and $u\ge0$ in $B_2$ (the ball with radius 2 centered at the origin) then
\begin{equation}\label{whi}
\inf_{B_{1}} u\ge c{\left(\int_{B_{1}}u^\varepsilon\right)}^{{1}/{\varepsilon}},
\end{equation}
where $ \ep(n, \lambda, \Lambda)>0$ and $c(n, \lambda, \Lambda,  q, \|b\|_{L^q(B_2)})>0$.
The weak Harnack inequality can also be viewed as a quantification of the strong maximum principle, in the sense that if a nonnegative supersolution is bounded below by a constant $a>0$ on a subset of positive measure $\mu>0$ then it is bounded below everywhere by $ka$, where $k>0$ depends  on $n,\lambda, \Lambda,q, \|b\|_{L^q(B_2)}$, and $\mu$.\smallskip

Another way to quantify SMP is to assume that $-\Lmin[u]$ rather than $u$ is bounded below by a constant $a>0$ on a subset of positive measure $\mu>0$. Then again $u$ is bounded below everywhere by $ka$, where $k>0$ depends only on $n, \lambda, \Lambda,q,\|b\|_{L^q(B_2)}$, and $\mu$. Equivalently, $\Lmin u\le0$ and $u\ge0$ in $B_2$ imply
\begin{equation}\label{qsmp}
\inf_{B_{1}} u\ge c{\left(\int_{B_{1}}(-\Lmin u)^\varepsilon\right)}^{{1}/{\varepsilon}},
\end{equation}
where $\ep, c>0$ are positive constants which depend on $n, \lambda, \Lambda,  q, \|b\|_{L^q(B_2)}$.

We refer to \re{qsmp} as the quantitative strong maximum principle (QSMP). It was essentially proved, for strong supersolutions of linear equations with bounded coefficients, in Krylov's book \cite{Kbook} (see also \cite{K}). We do not know of a reference for equations with unbounded coefficients, although the result is probably known to the experts (a proof will be included below). This not so widely known quantification of the SMP was used for instance in the work of Berestycki-Nirenberg-Varadhan \cite{BNV}, where it played an important role in the study of first eigenvalues of elliptic operators in nondivergence form. A more precise quantitative version for divergence-form operators appeared  in Brezis-Cabr\'e \cite{BC}, in the study of some nonlinear equations without solutions. The QSMP was used also in \cite{AS}, as a basic tool in the method developed there for proving nonexistence theorems for nonlinear elliptic inequalities.

The  interior estimates \eqref{whi} and \eqref{qsmp} have extensions to nonhomogeneous inequalities, $\Lmin[u]\le f(x)$, in which the right-hand sides of \eqref{whi} and \eqref{qsmp} have to be corrected by substracting $C\|f^+\|_{L^p(B_2)}$. As usual, we denote with $f^+$ (resp. $f^-$) the positive (resp. negative) part of $f$. \smallskip

It is our goal here to show that these interior estimates have boundary extensions. Theorems \ref{bqsmp}-\ref{bwhi} below quantify the Hopf lemma exactly like the WHI and QSMP quantify the SMP.

 We denote with $B_R^+=\{x\in\rn\::\:|x|<R,\;x_n>0\}$ a half-ball with a flat portion of the boundary included in $\{x_n=0\}$. We set $B_R^0=\partial B_R^+\cap\{x_n=0\}$ and $\lpr(B_R^+):= L^p(B_R^+)\cap L^q(B_R^+\cap\{x_n<d_0\})$, for some fixed $d_0>0$.
All the results can be stated for a bounded $C^{1,1}$-domain, see Corollary \ref{Bqsmpcor} and Theorem \ref{BWHI}.

In all that follows viscosity means $L^p$-viscosity in the sense of \cite{CCKS} -- see the next section for the definition and main properties of this notion. The theorems below are valid, with the same proofs, for $C$-viscosity sub- or super- solutions (as defined in \cite{CIL}), provided $b$ and $f$ are assumed continuous.

\begin{thm}\label{bqsmp} (boundary quantitative SMP, BQSMP) Assume that $u$ is a viscosity supersolution of $\Lmin[u]\le f$, $u\ge0$ in $B_2^+$, and $f\in \lpr(B_2^+)$. Then there exist constants $\varepsilon,c,C>0$ depending on $n$, $\lambda$, $\Lambda$, $p$, $q$, $d_0$ and $\|b\|_{L^q(B_2^+)}$, such that
 \begin{equation}\label{ineq1}
\inf_{B_1^+} \frac{u}{x_n} \ge c\left( \int_{B_{3/2}^+} (f^-)^{\varepsilon}\right)^{1/{\varepsilon}} - C\|f^+\|_{\lpr(B_2^+)}.
\end{equation}
\end{thm}

\begin{thm}\label{bwhi} (boundary weak Harnack inequality, BWHI) Assume that $u$ is a viscosity supersolution of $\Lmin[u]\le f$, $u\ge0$ in $B_2^+$, and $f\in \lpr(B_2^+)$. Then there exist constants $\varepsilon,c,C>0$ depending on $n$, $\lambda$, $\Lambda$, $p$, $q$, $d_0$ and $\|b\|_{L^q(B_2^+)}$, such that
 \begin{equation}\label{ineq2}
\inf_{B_1^+} \frac{u}{x_n} \ge c\left( \int_{B_{3/2}^+} \left(\frac{u}{x_n}\right)^\varepsilon\right)^{1/\varepsilon}-C \|f^+\|_{\lpr(B_2^+)}.
\end{equation}
\end{thm}

We also record the following  simple boundary extension of the local maximum principle for subsolutions (see \cite[Theorem 9.20]{GT}, \cite{CC}, \cite{KS3}).

\begin{thm}\label{blmp} (boundary local maximum principle, BLMP) Assume that $u$ is a viscosity subsolution of $\Lpl[u]\ge - f$ in $B_2^+$, $u\le0$ on $B_2^0$, and $f\in \lpr(B_2^+)$.  Then for each $r>0$ there exists $C>0$ depending on  $n$, $\lambda$, $\Lambda$, $p$, $q$, $r$, $d_0$, $\|b\|_{L^q(B_2^+)}$, such that
 \begin{equation}\label{ineq3}
\sup_{B_{1}^+} \frac{u^+}{x_n}\le C\left(\left( \int_{B_{3/2}^+} (u^+)^r\right)^{1/r} + \|f^+\|_{\lpr(B_2^+)}\right).
\end{equation}
\end{thm}

It is important to note that the combination of BWHI and BLMP yields a {\it boundary Harnack inequality for non-homogeneous equations}.

\begin{thm}\label{bhi} (inhomogeneous Harnack inequality, IHI) Assume that $u$ is a viscosity solution of $\Lmin[u]\le f^+$, $\Lpl[u]\ge -f^-$ in $B_2^+$, $u> 0$ in $B_2^+$, $u=0$ on $B_2^0$, and $f\in \lpr(B_2^+)$. Then there exists $C>0$ depending on $n$, $\lambda$, $\Lambda$, $p$, $q$, $d_0$ and  $\|b\|_{L^q(B_2^+)}$, such that
 \begin{equation}\label{fullbh}
\sup_{B_{1}^+} \frac{u}{x_n}\le C \left( \inf_{B_{3/2}^+} \frac{u}{x_n} + \|f\|_{\lpr(B_2^+)} \right).
\end{equation}
\end{thm}

When $f=0$ this is a  fundamental and very well known result which goes back to the work of Bauman \cite{B}, and the earlier work \cite{CCCC} for divergence form equations (see  Theorem \ref{safbh} below, and the references preceeding it). However, in all texts where this inequality  appeared it was proved by a method different from the above splitting into separate results for supersolutions and subsolutions (while the usual presentations of the {\it interior} full Harnack estimate include such splitting, see \cite{GT}, \cite{CC}). Furthermore, we do not actually know of any previous reference for \eqref{fullbh} with $f\not=0$, for equations in non-divergence form.

Theorems \ref{bqsmp} and \ref{bwhi} were first stated without proof in \cite{S2}, where they were used as a tool for establishing uniform a priori bounds for systems of elliptic inequalities.  I was recently informed by D. Moreira that other applications of Theorem \ref{bwhi}  to regularity theory and free boundary problems are in preparation \cite{BEM}, \cite{BMW}. The latter preprint provides an alternative and simpler proof of a weaker version of the BWHI (with $u/x_n$ replaced by $u$ and an interior integral in the right-hand side of \re{ineq2}).

To our knowledge, the only boundary result in the vein of Theorem \ref{bwhi} which appeared prior to this work is due to Caffarelli, Li and Nirenberg in the important work \cite{CLN} (Lemma 1.6 in that paper). With a  different proof (specific to linear equations with bounded coefficients), that result can be viewed as a weaker form of the ``growth" Lemma \ref{deggrowth} below, i.e. a first step to the proof of \re{ineq2} with $u/x_n$ replaced by $u$ in the integral.

We do not know of previous boundary estimates similar to Theorem \ref{bqsmp}.\smallskip

It is also worth noting that simpler but important inequalities also called boundary weak Harnack inequalities have appeared in the literature. In those inequalities $\inf u$ is bounded below by the integral of  $u_m^- = \min\{u,m\}$ where $m=\inf_{x_n=0}u$. That formulation is sufficient for a proof of H\"older regularity up to the boundary (see for instance \cite[Theorems 8.26 and 9.27]{GT}, \cite{KS2}, \cite{ARV}, \cite{SirARMA}), but is obviously void for a function which vanishes on the boundary ($m=0$) and does not imply the Hopf lemma or the full boundary Harnack inequality.

The proofs of Theorems \ref{bqsmp} and \ref{bwhi}  combine the proofs of the interior estimates with tools for boundary estimates such as Lipschitz bounds and global $W^{2,q}$-estimates, Hopf lemma and boundary barriers. In the following we will give complete and self-contained proofs,  both in order to provide a full quotable source for the boundary estimates which includes the interior estimates as particular cases, and because of the expected wide use of the theorems above. We will also save the interested reader the rather deep and somewhat hard to read treatment of the interior QSMP in the book \cite{Kbook}.

In the end, we observe that the optimal values of $\varepsilon$ in the WHI and QSMP, even for the Pucci operators, are an important open question. It is known that for {\it divergence} form operators, such as the Laplacian,  the interior WHI holds for any $\varepsilon<n/(n-2)$ (see \cite{GT}) and the interior QSMP holds for $\varepsilon =1$ (see \cite{BC}).\medskip

\noindent {\it Acknowledgement}. The author thanks warmly the anonymous referees, whose numerous remarks led to a substantial improvement of the presentation.

\section{Preliminaries}

We begin by recalling the notion of $L^p$-viscosity solution from \cite{CCKS}, adapted to equations with discontinuous coefficients. We recall we assume $p>p_0>n/2$. \smallskip

\noindent {\bf Definition 1.} Let $\Omega$ be a domain, $b,f\in L^p_{\mathrm{loc}}(\Omega)$. A function $u\in C(\Omega)$ is a  $L^p$-viscosity subsolution (resp. supersolution) of
$$
\Lpm[u]=f(x),
$$
we sometimes write $\Lpm[u]\ge(\mathrm{resp.}\le)f(x)$, provided
$$
\mathrm{ess}\liminf_{y\to x_0} (\Lpm[\phi(y)]-f(y))\le0
\quad
(\mathrm{resp.}\quad \mathrm{ess}\limsup_{y\to x_0} (\Lpm[\phi(y)]-f(y)\ge0),
$$
whenever $\phi$ is a function in $W^{2,p}_{\mathrm{loc}}(\Omega)$ such that $u-\phi$ attains local maximum (resp. minimum) at some point $x_0\in \Omega$. A solution is a function which is both a subsolution and a supersolution.\medskip

\noindent {\bf Definition 2.} Let $\Omega$ be a domain, $b,f\in L^p_{\mathrm{loc}}(\Omega)$. A function $u\in W^{2,p}_{\mathrm{loc}}(\Omega)$ is a  $L^p$-strong subsolution (resp. supersolution) of
$
\Lpm[u]=f(x),
$
provided
$$
\Lpm[u]\ge f(x)\quad \mbox{ a.e.  in }\; \Omega
\qquad
(\mathrm{resp.}\quad \Lpm[u]\le f(x)\quad \mbox{ a.e.  in }\; \Omega).
$$

We recall that strong (sub-)solutions are viscosity (sub-)solutions, and  a viscosity (sub-)solution which belongs to $ W^{2,p}_{\mathrm{loc}}(\Omega)$ is a strong (sub-)solution. We also recall that the comparison principle holds for $\mathcal{L}^\pm$, in the sense that if a subsolution is below a supersolution  on the boundary of a bounded $\Omega$, and one of them is strong, then they compare in the whole domain $\Omega$. For details, see \cite{CCKS}.\smallskip

In the rest of this section we quote a number of elliptic estimates which will be used in the proofs of Theorems \ref{bqsmp} and \ref{bwhi}. In all that follows $\Omega$ is a bounded $C^{1,1}$-domain.\smallskip

We begin with the ABP estimate, see \cite[Theorem 9.1]{GT}, in a version proved by Cabr\'e \cite{Cab}. The following is a particular case of Theorem 8.1 in \cite{KS2}, where inequalities with unbounded coefficients are considered.
\begin{thm}\label{ksabp} Assume $f\in L^p(\Omega)$, $q>n$, $q\ge p>p_0$, and $u$ is a viscosity solution of  $\Lpl[u]\ge -f$ in~$\Omega$. If for some $R_0,\sigma_0>0$ we have for each $x\in \Omega$
$$
|B_{R_0}(x)\setminus\Omega|\ge \sigma_0|B_{R_0}(x)|,
$$
then
$$
\sup_\Omega u \le \sup_{\partial \Omega}u^+ + CR_0^{2-n/p}\|f^+\|_{L^p(\Omega)},
$$
where $C=C(n,\lambda,\Lambda, p,q, \|b\|_{L^q(\Omega)},\sigma_0)$.
\end{thm}

We quote the following existence result and $W^{2,q}$-bound for equations with coefficients in $L^q$, included for instance in Proposition 2.4 in \cite{KS1} (see also \cite{W}).
\begin{thm}\label{ksexist} If $f\in L^q(\Omega)$, $\psi\in W^{2,q}(\Omega)$, $q>n$, then the equation $${\mathcal{L}}^\pm[u]= f\;\mbox{ in }\;\Omega, \qquad u=\psi\;\mbox{ on }\;\partial \Omega,$$ has a unique (among all viscosity solutions) solution which is in $W^{2,q}(\Omega)$ and
$$
\|u\|_{W^{2,q}(\Omega)}\le C (\|\psi\|_{W^{2,q}(\Omega)}+\|f\|_{L^q(\Omega)}),
$$
where $C=C(n, \lambda,\Lambda, q, \|b\|_{L^q(\Omega)},\Omega)$.
\end{thm}

We can easily  infer the following Lipschitz bound for inequalities of our type. We recall we set $\Omega_{d_0}:=\{x\in\Omega\,:\, \mathrm{dist}(x,\partial\Omega)<d_0\}$ where $d_0$ is a fixed small number, $\lpr(\Omega):= L^p(\Omega)\cap L^q(\Omega_{d_0})$, and (with a slight abuse of notation) $\lpr(B_R^+):= L^p(B_R^+)\cap L^q(B_R^+\cap\{x_n<d_0\})$.

\begin{thm}\label{lip}  1. Let $f\in \lpr(\Omega)$ and $u$ be a viscosity solution of  $\Lpl[u]\ge -f$ in~$\Omega$, $u\le0$ on $\partial \Omega$. Then
$$
u(x)\le C\|f^+\|_{\lpr(\Omega)}\, \mathrm{dist}(x,\partial \Omega),\quad x\in\Omega,
$$
 for some constant $C$ which depends on $n$, $\lambda$, $\Lambda$, $p$, $q$, $d_0$, $\|b\|_{L^q(\Omega)}$, and $\Omega$. \\ 2. If $f\in \lpr(B_2^+)$, and $u$ is a viscosity solution of $\Lpl[u]\ge -f$ in $B_2^+$, with $u\le0$ on $B_2^0$, then
$$
u\le C(\sup_{B_{3/2}^+}u^++\|f^+\|_{\lpr(B_2^+)})\, x_n \qquad \mbox{in }\;B_1^+.
$$
\end{thm}

\noindent {\it Proof.}  By the ABP inequality (Theorem \ref{ksabp} above) we know that in $\Omega$
$$
u\le C_0\|f^+\|_{L^p(\Omega)}.
$$
Let $h$ be a smooth function in $\Omega$ such that $h= C_0\|f^+\|_{L^p(\Omega)}$ in $\Omega\setminus \Omega_{d_0}$ and $h=0$ on $\partial \Omega$, with $\|h\|_{C^2(\Omega)}\le C(d_0) C_0\|f^+\|_{L^p(\Omega)}$. Let $v$ be the strong solution of $\Lpl[v] = -f^+ $ in $\Omega_{d_0}$, with $v=h$ on $\partial \Omega_{d_0}$, given by Theorem \ref{ksexist}. Then by the comparison principle $u\le v$ in $\Omega_{d_0}$, while by Theorem \ref{ksexist}
$$
\|v\|_{C^1(\Omega_{d_0/2})} \le C \|v\|_{W^{2,q}(\Omega_{d_0})}\le C (\|h\|_{W^{2,q}(\Omega_{d_0})} + \|f^+\|_{L^q(\Omega_{d_0})}) \le C\|f^+\|_{\lpr(\Omega)}.
$$
The proof of the second statement is similar.\hfill
$\Box$\smallskip

The next result is the interior weak Harnack inequality for viscosity supersolutions, in its version proved in \cite{KS2}. This inequality obviously implies the SMP.

\begin{thm}\label{kswh} If $f\in L^p(B_2)$ and $u$ is a viscosity solution of $\Lmin[u]\le f$, $u\ge0$ in $B_2$, then for some  $\ep= \ep(n, \lambda, \Lambda)>0$ and $c,C>0$ depending on $n, \lambda, \Lambda, p,q, \|b\|_{L^q(B_2)}$, we have
$$
\inf_{B_{1}} u\ge c{\left(\int_{B_{1}}u^\varepsilon\right)}^{{1}/{\varepsilon}}- C\|f^+\|_{L^p(B_2)}.
$$
\end{thm}

We also recall the local maximum principle for viscosity subsolutions, \cite[Theorem 9.20]{GT}, \cite{CC}, \cite{KS3}.

\begin{thm}\label{kslocmax} If $f\in L^p(B_2)$ and  $u$ is a viscosity solution of $\Lpl[u]\ge - f$ in~$B_2$, then for every $r>0$ there exists  $C>0$ depending on $n, \lambda, \Lambda,  p,q,r$, $\|b\|_{L^q(B_2)}$, such that
$$
\sup_{B_{1}} u^+\le C{\left(\int_{B_{3/2}}(u^+)^r\right)}^{{1}/{r}}+ C\|f^+\|_{L^p(B_2)}.
$$
\end{thm}

We now give the (simple) proof of Theorem \ref{blmp}.\smallskip

\noindent{\it Proof of Theorem \ref{blmp}}. Since the maximum of subsolutions is a viscosity subsolution the function $u^+=\max\{u,0\}$ solves $-\Lpl[u^+]\le f^+$ in $B_2^+$. Then extending $u^+$ and $f$  as zero in $B_2\setminus B_2^+$ we get a subsolution in $B_2$, to which we apply first Theorem \ref{kslocmax} and then Theorem~\ref{lip}. Observe that the radii $1$, $3/2$, $2$ in these theorems can be replaced by any other radii $R_1<R_2<R_3$, by rescaling.\hfill \Box\bigskip

Finally, we recall the following Carleson and boundary Harnack inequalities for {\it solutions} of homogeneous elliptic equations, which was first proved in \cite{CCCC}, \cite{B}. See also \cite[Proposition A.2]{ASS}, and, for equations with unbounded coefficients,  \cite{Saf}.

\begin{thm}\label{safbh} If $u$ is a strong solution of $\Lmin[u]\le 0\le \Lpl[u]$, $u>0$ in $B_2^+$,  $u=0$ on $B_2^0$ then for some $c_1,c_2,c_3>0$ depending on $\lambda, \Lambda, n, p,q, \|b\|_{L^q(B_2^+)}$, we have
$$
\inf_{B_1^+} \frac{u}{x_n} \ge c_1 u(0,\ldots,0,1/2)\ge c_2 \sup_{B_1^+} u,
$$
and if $v$ is another strong solution of $\Lmin[v]\le 0\le \Lpl[v]$, $v>0$ in $B_2^+$, and $v=0$ on $B_2^0$ then
$$
\inf_{B_1^+} \frac{u}{v} \ge c_3 \sup_{B_1^+} \frac{u}{v}.
$$
\end{thm}

This theorem easily implies the {\it Hopf lemma} for viscosity solutions of our inequalities.

\begin{thm}\label{hopfus} If $v$ is a viscosity supersolution of $\Lmin[v]\le 0$, $v>0$ in $B_2^+$,  then
$$
\inf_{B_1^+} \frac{v}{x_n} >0.
$$
\end{thm}
\smallskip

\noindent{\it Proof.}   Let $u$ be the strong solution of $\Lmin[u]=0$ in $B_{3/2}^+$, $u=v$ on $\partial B_{3/2}^+$ (given for instance by Theorem 7.1 in \cite{KS2}). By the comparison and the strong maximum principle we clearly have $v\ge u>0$ in $B_{3/2}^+$. So by Theorem~\ref{safbh} (with $B_2^+$ replaced by $B_{3/2}^+$) $\inf_{B_1^+} {v}/{x_n}\ge \inf_{B_1^+} {u}/{x_n}\ge c_2 \sup_{B_1^+} u >0$.\hfill $\Box$
\section{Proof of Theorem \ref{bqsmp}}

We start by observing it is sufficient to prove the following result. In the following the dependence of the constants in $\Omega$ is through its diameter and an upper bound of the curvature of $\partial\Omega$.
\begin{thm} \label{Bqsmp} There exist $\varepsilon,c>0$ depending on $n$,  $\lambda$, $\Lambda$, $q$, $\|b\|_{L^q}$ and $\Omega$ such that for each strong solution of
$\Lmin[u]\le 0$, $u\ge0$ in $\Omega$ we have

$$
\inf_{\Omega} \frac{u}{d} \ge c\left( \int_{\Omega} (-\Lmin[u])^{\varepsilon}\right)^{1/{\varepsilon}},
$$
where $d(x)=\mathrm{dist}(x,\partial \Omega)$.
\end{thm}

Indeed, this theorem and the boundary Lipschitz bound easily imply the following corollary.

\begin{cor} \label{Bqsmpcor} There exist $\varepsilon,c, C>0$ depending on $n$,  $\lambda$, $\Lambda$, $q$, $d_0$, $\|b\|_{L^q(\Omega)}$ and $\Omega$ such that for each viscosity solution of
$\Lmin[u]\le f\in \lpr(\Omega)$, $u\ge0$ in~$\Omega$, we have
$$
\inf_{\Omega} \frac{u}{d} \ge c\left( \int_{\Omega} (f^-)^{\varepsilon}\right)^{1/{\varepsilon}} - C\|f^+\|_{\lpr(\Omega)}.
$$
\end{cor}

\noindent {\it Proof.} We take the (nonnegative) strong solutions $v,w$ of $\Lmin[v]= - f^-$ in $\Omega$, $v=u$ on $\partial\Omega$, and $\Lpl[w]= - f^+$ in $\Omega$, $w=0$ on $\partial\Omega$. Then $$\Lmin[v] \ge \Lmin[u] + \Lpl[w] \ge \Lmin[u+w]\quad\mbox{ in }\Omega, $$
 so  by the comparison principle $v\le u+w$ in $\Omega$. Hence, by applying Theorem~\ref{Bqsmp} to $v$,
 $$
c\left( \int_{\Omega} (f^-)^{\varepsilon}\right)^{1/{\varepsilon}} \le \inf_{\Omega} \frac{v}{d}\le \inf_{\Omega} \frac{u}{d} + \sup_{\Omega} \frac{w}{d}\le \inf_{\Omega} \frac{u}{d} + C \|f^+\|_{\lpr(\Omega)},
$$
where we used the Lipschitz bound, Theorem \ref{lip} applied to $w$. \hfill $\Box$\bigskip

\noindent {\it Proof of Theorem \ref{bqsmp}.}
Apply the previous corollary in some smooth domain $\Omega$ such that $B_{3/2}^+\subset\Omega\subset B_2^+$. \hfill $\Box$\bigskip

In the rest of this section we prove Theorem \ref{Bqsmp}. The overall scheme of the proof is similar to the one used by Krylov in the proof of the interior estimate, but we need to provide boundary extensions to all steps of that proof. For the reader's convenience we give a complete proof of Theorem~\ref{Bqsmp}, which encompasses both the interior and the boundary estimate. \bigskip

\noindent {\it Proof of Theorem \ref{Bqsmp}}. To simplify notations, we will assume $\Omega = B_1$ and $d=d(x)= \mathrm{dist}(x,\partial B_1)$. For general $\Omega$ we can either use coverings by balls or repeat the arguments below replacing $B_1$ by $\Omega$, with obvious changes.

 By the SMP and the Hopf lemma (Theorems \ref{kswh} and \ref{hopfus}) we have either $u\equiv 0$ or $\inf_{B_{1}} u/d>0$. Dividing $u$ by $\inf_{B_{1}} u/d$ we can assume $\inf_{B_{1}} u/d=1$. We are going to show that there exist positive constants $\ep_0,C_0$ depending only on $n$,  $\lambda$, $\Lambda$, $q$, and $\|b\|_{L^q(B_1)}$, such that
\begin{equation}\label{meas1}
|\{-\Lmin[u] \ge t\}\cap B_1|\le C_0\min\{1,t^{-2\ep_0}\}.
\end{equation}
Then the theorem follows, since we can write
\begin{eqnarray*}
\int_{B_1} (-\Lmin[u])^{\ep_0} &=& \ep_0 \int_0^\infty t^{\ep_0-1}|\{-\Lmin[u] \ge t\}\cap B_1|\,dt\\
 &\le& C_0 \ep_0\int_0^\infty t^{\ep_0-1} \min\{1,t^{-2\ep_0}\} =C.
\end{eqnarray*}

Next, observe that \re{meas1} is equivalent to the following statement: {\it there exist positive constants $A_0,c_0$ depending only on $n$,  $\lambda$, $\Lambda$, $q$,  and $\|b\|_{L^q}$, such that for each strong solution of
$\Lmin[u]\le 0$, $u\ge0$ in $B_1$,
\begin{equation}\label{meas2}
\mbox{if }\;\; s=|\{-\Lmin[u] \ge 1\}\cap B_1|/|B_1|\qquad\mbox{then}\qquad  \inf_{B_{1}} \frac{u}{d} \ge c_0
s^{A_0}.
\end{equation}}
We will prove this claim by an iteration procedure.

We first  record the following fact.
\begin{lem}\label{lemkr4} The solution of the problem
$$
\left\{
\begin{array}{rclcc}
\Lmin[v]&=&-1&\mbox{in}&B_1\\
v&=&0&\mbox{on}&\partial B_1,
\end{array}
\right.
$$
is such that
$$
v\ge c_0\, d $$
where $d(x)= \mathrm{dist}(x,\partial B_1)$ and $c_0>0$ depends on $n,\lambda, \Lambda, q$, and $\|b\|_{L^q(B_1)}$.
\end{lem}

We postpone the proof of this lemma.
\bigskip

The following result will permit to us to start the iteration.
\begin{lem}\label{lemkr1} There exist positive constants $\delta_1,c_1$ depending on $n$,  $\lambda$, $\Lambda$, $q$, $d_0$, and $\|b\|_{L^q}$, such that for each strong solution of
$\Lmin[u]\le 0$, $u\ge0$ in $B_1$
$$
 \frac{|\{-\Lmin[u] \ge 1\}\cap B_1|}{|B_1|}\ge 1-\delta_1\qquad\mbox{implies}\qquad  \inf_{B_{1}} \frac{u}{d} \ge c_1.
$$
\end{lem}

\noindent{\it Proof.} Let $A=\{-\Lmin[u] \ge 1\}\cap B_1$ and let $v$ and $w$ be the strong solutions of the Dirichlet problems
$$
\left\{
\begin{array}{rclcc}
\Lmin[v]&=&-1&\mbox{in}&B_1\\
v&=&0&\mbox{on}&\partial B_1,
\end{array}
\right.
\qquad
\left\{
\begin{array}{rclcc}
\Lpl[w]&=&-\chi_{B_1\setminus A}&\mbox{in}&B_1\\
w&=&0&\mbox{on}&\partial B_1,
\end{array}
\right.
$$
given by Theorem \ref{ksexist} (here and in the following $\chi_Z$ denotes the characteristic function of a set $Z$). Then
$$
\Lmin[u]\le \Lmin[v] - \Lpl[w] \le \Lmin[v-w]
$$
in $B_1$ so by the comparison principle $ u\ge v-w$ in $B_1$. By the previous lemma $v\ge c_0 d$, while the $C^{1}$-bound given by Theorem \ref{ksexist} (or the Lipschitz bound, Theorem \ref{lip}),  implies
$$
w\le C_0|B_1\setminus A|^{1/q}d\le C_0\delta_1^{1/q}d.
$$
The lemma follows, with $\delta_1=[c_0/(2C_0)]^{q}$. \hfill $\Box$
\bigskip

The following lemma was essentially proved by  Krylov in his treatment of the interior estimate in \cite{Kbook}. We give the proof for completeness.

\begin{lem}\label{lemkr2}
Let $G$ be an open set in $B_1$, and $f,g_in L^p(G)$ be nonnegative functions on $B_1$, with $f,g\not\equiv0$ on $B_1$ and  $g\equiv0$ on $B_1\setminus G$.

Suppose that for each $x_0\in G$ there exists a ball $B_\rho(x_1)\subset B_1$ such that $x_0\in B_\rho(x_1)$ and  the (unique strong) solutions of the Dirichlet problems
$$
\left\{
\begin{array}{rclcc}
\Lmin[u_1]&=&-f&\mbox{in}&B_\rho(x_1)\\
u_1&=&0&\mbox{on}&\partial B_\rho(x_1),
\end{array}
\right.
\qquad
\left\{
\begin{array}{rclcc}
\Lpl[v_1]&=&-g&\mbox{in}&B_\rho(x_1)\\
v_1&=&0&\mbox{on}&\partial B_\rho(x_1),
\end{array}
\right.
$$
verify $u_1(x_0)\ge v_1(x_0)$.

Then the solutions  of the Dirichlet problems
$$
\left\{
\begin{array}{rclcc}
\Lmin[u]&=&-f&\mbox{in}&B_1\\
u&=&0&\mbox{on}&\partial B_1,
\end{array}
\right.
\qquad
\left\{
\begin{array}{rclcc}
\Lmin[v]&=&-g&\mbox{in}&B_1\\
v&=&0&\mbox{on}&\partial B_1,
\end{array}
\right.
$$
are such that $u\ge v$ in $B_1$.
 In other words, increase on small scales implies increase on the whole domain.
\end{lem}

\noindent{\it Proof.} Assume the lemma is proved for any function $g$ which vanishes outside some closed set contained in $G$. Then the lemma follows, since the solutions of the Dirichlet problem with  right-hand side $g\chi_{A_k}$ converge uniformly to the solution of the problem with right-hand side $g$, by the ABP inequality (here $A_k\nearrow G$ is a sequence of closed sets).

So we can assume $g$ vanishes outside some closed $\Gamma\subset G$. Let $\delta>0$ be arbitrary and $w:=v-(1+\delta)u$. Assume for contradiction that $$\alpha:= \max_{B_1}w>0.$$

We claim that this maximum is attained at some point $x_0\in \Gamma$. Indeed, if it is attained in the domain $B_1\setminus\Gamma$ we can apply the SMP to the inequality $$\Lmin[\alpha-w]=-\Lpl[w]\le (1+\delta)\Lmin[u]-\Lmin[v]  \le 0,$$ which holds in $B_1\setminus\Gamma$, and deduce that $w\equiv \alpha$ in some connected component of $B_1\setminus\Gamma$. Since the boundary of this connected component may contain only points on $\partial \Omega$ (where $w=0$) and points on $\partial \Gamma$, we conclude that $w=\alpha$ for some point on $\partial \Gamma\subset \Gamma$.

For this point $x_0$ take the ball $B_\rho(x_1)$ and the functions $u_1$, $v_1$ given by the assumption of the lemma. Then in $B_\rho(x_1)$
$$
\Lmin[v-v_1]\ge \Lmin[v]-\Lpl[v_1]=0= \Lmin[u]-\Lmin[u_1]\ge \Lmin[u-u_1]
$$
while on $\partial B_\rho(x_1)$
$$
v-v_1=v\le \alpha + (1+\delta)u  =\alpha + (1+\delta)(u-u_1)
$$
so by the maximum principle
$$
v(x_0)-v_1(x_0)\le \alpha + (1+\delta) (u(x_0)-u_1(x_0))
$$
that is, $(1+\delta)u_1(x_0)\le v_1(x_0)$, a contradiction with the assumption of the lemma.

We have proved that $(1+\delta)u\ge v$ in $B_1$ for each $\delta>0$, and we conclude by letting $\delta\to0$.\hfill $\Box$
\bigskip

We also quote the following well known measure theoretic result, Krylov's ``propagating ink spots lemma'' (see for instance Lemma 1.1  in \cite{Saf80}, or Lemma~6 on page 122 of \cite{Kbook}).
\begin{lem}\label{lemkr3}
Let  $A\subset B_1$ be such that $|A| \leq \eta |B_1|$ where $\eta \in (0,1)$. Consider the family of balls
$$
\mathcal{F} = \{B_\rho(x)  \; ball\, , B_\rho(x)\subset B_1\,, \,|B_\rho(x)\cap A|\geq \eta |B_\rho|\, \}.
$$
Then there exist $\xi = \xi(n,\eta)<1$ and $\zeta= \zeta(n,\eta)>0$ such that
$$|\tilde{A}|=\left|\bigcup_{B_\rho(x)\in\mathcal{F}} B_{\xi\rho}(x)\right| \geq (1+\zeta)|A|.$$
\end{lem}

We now give the proof of Theorem \ref{Bqsmp}. We recall it was reduced to \re{meas2}.\smallskip

\noindent{\it Proof of \re{meas2}}. Given a measurable subset $A\subset B_1$ with positive measure, we denote
\begin{equation}\label{eqiw}
r(A) = \inf_{B_{1}} \frac{w}{d}>0, \qquad \mbox{where}\qquad
\left\{
\begin{array}{rclcc}
\Lmin[w]&=&-\chi_A&\mbox{in}&B_1\\
w&=&0&\mbox{on}&\partial B_1
\end{array}
\right.
\end{equation}
($w$ is given by Theorem \ref{ksexist}; $r(A)>0$ by the SMP and the Hopf lemma, Theorems \ref{kswh} and \ref{hopfus}), and
$$
\mu(s) := \inf_{A:\, |A|\geq s|B_1|} r(A).
$$

We observe that   $r(A)\le C_1$ for some constant $C_1$ depending only on $n,\lambda,\Lambda, p,q$, $\|b\|_{L^q(B_1)}$, by the Lipschitz bound (Theorem \ref{lip}).

  We need to prove that there exist $c_0,A_0>0$ depending on $n$, $\lambda$, $\Lambda$, $q$, and $\|b\|_{L^q(B_1)}$, for which
\begin{equation}\label{meas3}
\mu(s)\ge c_0s^{A_0}, \qquad \mbox{for }\;s\in(0,1).
\end{equation}
Indeed, if this is proved, by the comparison principle $u\ge w$ in $B_1$, where $w$ is the solution of \re{eqiw} with $A = \{-\Lmin[u]\ge1\}\cap B_1$, and we infer \re{meas2}.

Let us prove \re{meas3}. First,
Lemma \ref{lemkr1} implies that $\mu(s)\ge c_1\ge c_1 s>0$ provided $s\in (1-\delta_1,1]$. Second, it is enough to find $\zeta,\bar C>0$ depending only on $n$, $\lambda$, $\Lambda$, $q$,   $\|b\|_{L^q}$, such that for any $0<s_1<s_2\le1$
\begin{equation}\label{meas33}
\mu(s_2)>0\;\mbox{ and } \;s_1= (1+\zeta)^{-1}s_2\qquad\mbox{imply}\qquad \mu(s_1)\ge {\bar C}^{-1} \mu(s_2).
\end{equation}
Indeed, a simple iteration argument shows that the nondecreasing function $\mu(s)$ satisfies \re{meas3} provided it satisfies \re{meas33} -- since then for all $k\ge 1$ $$\mu\left((1+\zeta)^{-k}\right)\ge {\bar C}^{-k}\mu(1).$$

We set $\zeta>0$ to be the constant from Lemma \ref{lemkr3} applied with $\eta=1-\delta_1$, where $\delta_1$ is the constant from Lemma \ref{lemkr1}.

Let $0<s_1<s_2\le1$ be such that $s_1= (1+\zeta)^{-1}s_2$ and $\mu(s_2)>0$. In order to prove that $\mu(s_1)\ge {\bar C}^{-1} \mu(s_2)$, we need to show that for every  subset $A\subset B_1$ such that $|A|\ge s_1|B_1|$, there exists a subset $\tilde A\subset B_1$ such that $|\tilde{A}|\ge s_2|B_1|$ and $r(A)\ge {\bar C}^{-1} r(\tilde A)$ ($\bar C$ is still to be chosen).

Let $A\subset B_1$ be such that $|A|\ge s_1|B_1|$. First, if  $|A|\ge (1-\delta_1)|B_1|$ then by Lemma \ref{lemkr1} and the Lipschitz bound we get $r(A)\ge c_1\ge (c_1/C_1) r(\tilde A)$, for every $\tilde A\subset B_1$. Hence we can assume that  $|A|< (1-\delta_1)|B_1|$. Let now $\tilde{A}$ be the open set constructed in Lemma~\ref{lemkr3}, so that $|\tilde{A}|\ge (1+\zeta)|A|\ge (1+\zeta)s_1|B_1|= s_2|B_1|$.

Let  $u$ and $v$ be the solutions of the Dirichlet problems
$$
\left\{
\begin{array}{rclcc}
\Lmin[u]&=&-\bar{C}\chi_A&\mbox{in}&B_1\\
u&=&0&\mbox{on}&\partial B_1,
\end{array}
\right.
\qquad
\left\{
\begin{array}{rclcc}
\Lmin[v]&=&-\chi_{\tilde{A}}&\mbox{in}&B_1\\
v&=&0&\mbox{on}&\partial B_1.
\end{array}
\right.
$$
We claim that we can choose $\bar{C}$, depending only on $n$, $\lambda$, $\Lambda$, $q$,  and $\|b\|_{L^q}$, in such a way that $u\ge v$ in $B_1$ -- which implies that $r(A)\ge \bar C^{-1}r(\tilde A)$.

It remains to prove the last claim. For this we will use Lemma \ref{lemkr2} with $f= \bar C \chi_{A}$ and $g= \chi_{\tilde{A}}$. Let $x_0\in \tilde{A}$. Then there exists a ball $B_\rho(x_1)\subset B_1$ such that $x_0\in B_{\xi\rho}(x_1)$ for some $\xi<1$ given by Lemma \ref{lemkr3} and depending on the right quantities, and $$|A\cap B_\rho(x_1)|\ge (1-\delta_1)|B_\rho(x_1)|.$$
The latter inequality and Lemma \ref{lemkr1} imply that (after rescaling $x\to x/\rho$) the solution of the Dirichlet problem
$$
\left\{
\begin{array}{rclcc}
\Lmin[u_1]&=&-\chi_A&\mbox{in}&B_\rho(x_1)\\
u_1&=&0&\mbox{on}&\partial B_\rho(x_1),
\end{array}
\right.
$$
is such that
$$u_1\ge \tilde{c}\rho^2\quad\mbox{on}\quad B_{\xi\rho}(x_1),$$ for some $\tilde{c}>0$ which depends only on $n$, $q$, $\lambda$, $\Lambda$ and $\|b\|_{L^q}$ (recall $\xi$ only depends on these too).

On the other hand, the solution of
$$
\left\{
\begin{array}{rclcc}
\Lmin[v_1]&=&-\chi_{\tilde{A}}&\mbox{in}&B_\rho(x_1)\\
v_1&=&0&\mbox{on}&\partial B_\rho(x_1),
\end{array}
\right.
$$
is obviously such that
$$
v_1\le \tilde{C}\rho^2 \quad\mbox{on}\quad B_{\rho}(x_1),$$
by the ABP inequality.

Setting $\bar{C} = \tilde{C}/\tilde{c}$, the claim now follows from Lemma \ref{lemkr2}, and the proof is finished. \hfill $\Box$
\bigskip

We now return to Lemma \ref{lemkr4}.

\noindent{\it Proof of Lemma \ref{lemkr4}}. Assume there exists a sequence of functions $b_k$ such that $\|b_k\|_{L^q(B_1)}\le C$ and points $y_k\in B_1$ with $v_k(y_k)/d(y_k) \to 0$, where $v_k$ is the (strong) solution of
\begin{equation}\label{thi}
\left\{
\begin{array}{rclcc}
\mm(D^2v_k) - b_k|Dv_k|&=&-1&\mbox{in}&B_1\\
v_k&=&0&\mbox{on}&\partial B_1.
\end{array}
\right.
\end{equation}
 By the global $W^{2,q}$-estimate (Theorem \ref{ksexist}) a subsequence of $v_k$ converges to a function $v$ in $C^{1,\alpha}(B_1)$ for some $\alpha >0$.
\smallskip

\noindent {\it Step 1. We have the interior estimate: for each $r<1$ there exists $c_r>0$ such that $v\ge c_r$ in $B_r$.}

\noindent {\it Proof.} If a subsequence of  $y_k$ converges to a point $y_0\in B_1$, we apply the suitably rescaled interior weak Harnack inequality, Theorem \ref{kswh}, to the  Dirichlet problem \re{thi} in $B_R$ with $R<1$,
 $$
\int_{B_R} v_k^\varepsilon \le C_R \inf_{B_R} v_k^\varepsilon \le C_R (v_k(y_k))^\varepsilon ,
$$
pass to the limit as $k\to\infty$ for each fixed $R<1$, and deduce that $v\equiv0$ in $B_1$. But then $\mm(D^2v_k) =f_k$ where $f_k:=-1+b_k|Dv_k|\to -1$ in $L^q(B_1)$, so by the ABP inequality $v_k(x)$ converges uniformly to $(\lambda/2n)(1-|x|^2)$, a contradiction with $v_k(y_k)\to 0$.
\smallskip

\noindent {\it Step 2. We have the interior version of Theorem \ref{Bqsmp} : for each $r<1$ there exists $c_r>0$ such that for each strong solution of
$\Lmin[u]\le 0$, $u\ge0$ in $B_1$ we have

$$
\inf_{B_{r}} u \ge c_r\left( \int_{B_{1}} (-\Lmin[u])^{\varepsilon}\right)^{1/{\varepsilon}}.
$$}

\noindent {\it Proof.} By using only the Step 1 we just established, we can repeat almost verbatim the whole proof of Theorem \ref{Bqsmp}, provided we replace the conclusion of Lemma~\ref{lemkr1} by $\inf_{B_{r}} u \ge c_r$, the inequality $v\ge c_0 d$ by $v\ge c_r$ in $B_r$ in the proof of Lemma~\ref{lemkr1}, and define $r(A) = \inf_{B_r} w$.
\smallskip

\noindent {\it Step 3. Conclusion.}  Assume that a subsequence of $y_k$ converges to a point $y_0\in \partial B_1$ (say $y_0 = e =(0,\ldots,0,1))$, i.e. $\frac{\partial v_k}{\partial e}(e)\to 0$ as $k\to \infty$. We take
$w_k$ to be the (strong) solution of the problem
$$
\left\{
\begin{array}{rclcc}
\mm(D^2w_k) - b_k|Dw_k|&=&0&\mbox{in}&B_{1/2}(e/2)\\
w_k&=&z_k&\mbox{on}&\partial B_{1/2}(e/2),\\
\end{array}
\right.
$$
where $z_k$ is a smooth function on $B^\prime := B_{1/2}(e/2)$, such that $z_k=0$ on $\partial B^\prime\cap\{x\,:\,x_n>3/4\}$, $z_k=v_k$ on $\partial B^\prime\cap\{x\,:\,x_n<1/2\}$, $0\le z_k\le v_k$ on $\partial B^\prime$, and $\|z_k\|_{W^{2,q}(B^\prime)}\le C\|v_k\|_{W^{2,q}(B^\prime)}$.

By the comparison principle $w_k\le v_k$ in $B^\prime$. From Step 2 with $r=1/\sqrt{2}$ applied to $v_k$ we already know that $$w_k=v_k\ge c_0\;\mbox{ on }\;\partial B^\prime\cap\{x\,:\,x_n<1/2\}.$$
 By Theorem \ref{ksexist}
$$
\|w_k\|_{C^1(B^\prime)}\le C\|w_k\|_{W^{2,q}(B^\prime)} \le C\|v_k\|_{W^{2,q}(B^\prime)}\le C,
$$
hence there exists $d_0>0$ such that $$w_k\ge c_0/2\qquad\mbox{in }\; B^\prime\cap\{x\,:\,x_n<1/4\}\cap\{x\,:\,\mathrm{dist}(x,\partial B^\prime)<d_0\}.$$
By applying the interior weak Harnack inequality in $B_{\frac{1-d_0}{2}}(e/2)$ we infer that
 $w_k(e/2)\ge c_1>0$. Then by applying Theorem \ref{safbh} to $w_k$ (after straightening $\partial B^\prime\cap\{x\,:\,x_n>3/4\}$) we get the contradiction
$$
-\frac{\partial v_k}{\partial e}(e) \ge -\frac{\partial w_k}{\partial e}(e)\ge  c_2w_k(e/2) \ge c_2c_1>0.
$$
This proves Lemma \ref{lemkr4}. \hfill $\Box$

\section{Proof of Theorem \ref{bwhi}}

In this section we give the proof of  the boundary weak Harnack inequality, Theorem \ref{bwhi}.

In the following $Q_\rho=Q_\rho(\rho e)$ denotes the cube with center $\rho e$ and side $\rho$, where $e= (0,\ldots,0,1/2)$. To avoid confusion, the reader's attention is brought to the fact that $Q_\rho$ is not centered at the origin but has its bottom on $\{x_n=0\}$. We also recall $d_0$ is the width of a neighborhood of the lower boundary of the cube $Q_2$, in which $f\in L^q$. Without loss we assume $d_0\le1/4$.

\begin{thm} \label{BWHI} (BWHI) 1. There exist $\varepsilon,c,C>0$ depending on $n$, $p$, $q$, $\lambda$, $\Lambda$, $d_0$, and $\|b\|_{L^q(Q_2)}$ such that for each  viscosity solution of
$\Lmin[u]\le  f(x)$, $u\ge0$ in $Q_2$ we have
$$
\inf_{Q_{1}} \frac{u}{x_n} \ge c\left( \int_{Q_{1}} \left(\frac{u}{x_n}\right)^\varepsilon\right)^{1/{\varepsilon}} - C\|f^+\|_{\lpr(Q_2)}.
$$

2. There exist $\varepsilon,c,C>0$ depending on $n$, $p$, $q$, $\lambda$, $\Lambda$, $d_0$, $\|b\|_{L^q(\Omega)}$ and $\Omega$ such that for each  viscosity solution of
$\Lmin[u]\le  f(x)$, $u\ge0$ in $\Omega$ we have
$$
\inf_{\Omega} \frac{u}{d} \ge c\left( \int_{\Omega} \left(\frac{u}{d}\right)^\varepsilon\right)^{1/{\varepsilon}} - C\|f^+\|_{\lpr(\Omega)},
$$
where $\Omega$ is a  bounded $C^{1,1}$-domain and we set $d=d(x)=\mathrm{dist}(x,\partial \Omega)$.
\end{thm}

To prove this theorem, we first establish the following growth lemma.

\begin{lem}[growth lemma]\label{deggrowth} Given $\nu>0$, there exist $k,a>0$ depending on $n$, $p$, $q$, $\lambda$, $\Lambda$, $d_0$, $\nu$, $\|b\|_{L^q}$,  such that if $u$ is a viscosity solution of
$$\Lmin[u]\le  f(x), \quad f,u\ge0\mbox{ in }Q_2, \qquad\mbox{and}\qquad \|f\|_{\lpr(Q_2)} \le a,$$ and we have
$$
|\{u>x_n\}\cap Q_1|\ge \nu ,
$$
then $u> kx_n$ in $Q_1$.
\end{lem}

This lemma implies the weak boundary Harnack inequality, Theorem~\ref{BWHI}, through a cube decomposition procedure. Here are the details.\smallskip

\noindent {\it Proof of Theorem \ref{BWHI}}. The second part of the theorem is an easy consequence of the first, by locally straightening the boundary and covering it with balls in which such straightening is possible.

To prove the first part of the theorem we will show that there exist $M>0$, $\mu<1$ and $\delta_0>0$, depending on the appropriate quantities, such that if $u$ is a solution of $$\Lmin[u]\le  f(x),\quad f,u\ge0\quad\mbox{in}\; Q_2,\quad \|f\|_{\lpr(Q_2)}\le \delta_0,\quad \mbox{and}\quad \inf_{Q_1} (u/x_n) \le 1,$$ then
\begin{equation}\label{iter}
|\{u/x_n> M^j\}\cap Q_1|\le (1-\mu)^j.
\end{equation}
After the inequality \re{iter} is proved we infer from it that for some $\varepsilon_0>0$ we have  $|\{u/x_n> t\}\cap Q_1|\le C\min\{1,t^{-2\varepsilon_0}\}$ for $t>0$, and hence similarly to the proof of the BQSMP in the previous section (see \re{meas1})
$$\int_{Q_1} (u/x_n)^{\varepsilon_0}\le C.$$
For each $\beta>0$, we apply this inequality to the function $u$ from Theorem \ref{BWHI} divided by $\inf_{Q_1} u/x_n +\beta + \|f^+\|/\delta_0$, let $\beta\to0$ and infer the theorem.

We next prove \re{iter}. Let $0<k<1$ and $a>0$ be the numbers from Lemma \ref{deggrowth} for $$\nu=(d_0/4)^n/2.$$

 Let $u$ be a solution of $\Lmin[u]\le  f(x)$, $u\ge0$ in $Q_2$, $\|f\|_{\lpr(Q_2)}\le \delta_0$, such that $\inf_{Q_1} u/x_n \le 1$.

Set $M=1/k$. Replacing $u$ by $ku$ in Lemma \ref{deggrowth} (this does not make the norm of $f$ larger since $k<1$)  we see that $\inf_{Q_1} ku/x_n \le k$ implies
 $$|\{u/x_n>M\}\cap Q_1| = |\{ku> x_n\}\cap Q_1|  \le \nu<1/2,$$
so \eqref{iter} is true for $j=1$  for every $\mu\le 1/2$. We fix
$$
\mu=c_0/2,
$$
where  $c_0<1$ is the constant from the following (equivalent) version of the propagating ink spots lemma, Lemma \ref{lemkr3}.

\begin{lem}\label{eqink}
Let $A\subset B \subset Q_1$ be two measurable sets. Assume  there exists $\alpha>0$ such that $|A|\leq (1-\alpha) |Q_1|$, and for any $x_0\in Q_1$, $\rho>0$ such that the cube $ Q=Q_\rho(x_0)\subset Q_1$ we have
$$
\mbox{if }\;|Q\cap A|\geq (1-\alpha)|Q| \;\mbox{ then }\; \; Q\subset B.
$$
Then $$|A|\leq (1-c_0\alpha)|B|,
$$
for some constant $c_0=c_0(n)\in (0,1)$.
\end{lem}

We define the sets $A = \{u/x_n> M^j\}\cap Q_1$ and $B=\{u/x_n> M^{j-1}\}\cap Q_1$. Let us fix some cube $Q=Q_\rho(x_0)\subset Q_1$ ($x_0=(x_0^\prime,x_{0,n})$ is the center, $\rho$ is the side of this cube), such that $$|A\cap Q|\ge (1-\mu) |Q|.$$ We want to show that $Q\subset B$, i.e. $u/x_n> M^{j-1}$ in $Q$. Then Lemma \ref{eqink} and an induction argument easily imply \re{iter}.

We will distinguish several cases.\smallskip

{\it Case 1.}  $\rho\ge d_0/4$.  Then $$|A\cap Q_1|\ge|A\cap Q|\ge (1-\mu) |Q|\ge | Q|/2 = \rho^n/2 \ge \nu.$$ That is, $|\{u/M^j> x_n\}\cap Q_1|\ge \nu$; then by Lemma \ref{deggrowth} we get $u/x_n>M^{j-1}$ in $Q_1$ which contains $Q$. \smallskip

{\it Case 2.} $\rho<d_0/4$. We further divide the argument in two subcases.

{\it Case 2.1.} $\rho< d_0/4$ and $Q\cap\{x_n<d_0/2\}\not=\emptyset$.
 These hypotheses imply that  $f\in L^q$ in a neighborhood of $Q=Q_\rho(x_0)$.

We rescale the variables by setting $$y=(y^\prime, y_n) = T_{\rho^\prime}(x):= \frac{(x^\prime-x_0^\prime, x_n)}{\rho^\prime}\,,\qquad\mbox{where }\; \rho^\prime := 2x_{0,n}\in [\rho, 5d_0/4),
  $$
  and
  $$
  v(y) =\frac{ u(x)}{\rho^\prime}= \frac{1}{\rho^\prime}u(x_0^\prime + \rho^\prime y^\prime, \rho^\prime y_n),\qquad  \tilde b(y) = b(x), \quad \tilde f(y) = f(x).
  $$
  Then $v$ is a solution of
\begin{equation}\label{equs}
\mm(D^2v) - \rho^\prime \tilde b(y)|Dv| \le \rho^\prime \tilde f(y)
\end{equation}
in the cube $Q_{1/\rho^\prime}\supset Q_2$ (recall we assume $d_0\le1/4$, so $\rho^\prime<1/2$).

Observe that $T_{\rho^\prime}(Q)= Q_{\rho/\rho^\prime}(e)$, and
 $\|\rho^\prime \tilde f\|_{L^r(T_{\rho^\prime}(A))} = (\rho^\prime)^{1-n/r}\|f\|_{L^r(A)}$, for every $r\ge 1$ and every set $A$. So if $p\ge n$ we have $$\|\rho^\prime \tilde f\|_{L^\prime(Q_2)}\le\|\rho^\prime \tilde f\|_{L^\prime(Q_{1/\rho^\prime})}\le \|f\|_{L^\prime(Q_2)}\le \delta_0,$$
whereas if $p<n<q$
\begin{eqnarray*}
\|\rho^\prime \tilde f\|_{L^\prime(Q_2)} &=&
 \|\rho^\prime \tilde f\|_{L^q(\{y_n<d_0\}\cap Q_2)} + \|\rho^\prime \tilde f\|_{L^p(\{y_n>d_0\}\cap Q_2)} \\
 &\le& \|\rho^\prime \tilde f\|_{L^q(\{y_n<d_0\}\cap Q_2)} + |Q_2|^{1/p-1/q} \|\rho^\prime \tilde f\|_{L^q(\{d_0<y_n<d_0/\rho^\prime\}\cap Q_2)} \\
 & & + \|\rho^\prime \tilde f\|_{L^p(\{d_0/\rho^\prime<y_n<2\}\cap Q_2)}\\
 &\le& C (\rho^\prime)^{1-n/q} \|f\|_{L^q(\{x_n<d_0\}\cap Q_2)} +
(\rho^\prime)^{1-n/p} \|f\|_{L^p(\{d_0<x_n<2/\rho^\prime\}\cap Q_2)}\\
&\le& (C (5d_0/4)^{1-n/q} + (d_0/2)^{1-n/p}) \|f\|_{L^\prime(Q_2)} \le C \delta_0
\end{eqnarray*}
(the set in the last $L^p$-norm is empty if $\rho^\prime<d_0/2$). We set $\delta_0$ small enough so that $C \delta_0 <a$, where $a$ is the constant from Lemma \ref{deggrowth}.\smallskip

We have  $u(x)/x_n = v(y)/y_n$, so $|A\cap Q|\ge (1-\mu) |Q|$ is equivalent to
\begin{equation}\label{tukk}
|\{v/M^j > y_n\}\cap Q_{\rho/\rho^\prime}(e)| \ge (1-\mu)  |Q_{\rho/\rho^\prime}(e)|.
\end{equation}

 If $\rho\ge  \rho^\prime/4$ (i.e. $Q$ is ``close" to the bottom boundary) \re{tukk} implies
 $$
|\{v/M^j > y_n\}\cap Q_{1}| \ge (1-\mu)  4^{-n}  \ge (1/2) 4^{-n}  \ge \nu,
 $$
so  we can  apply Lemma \ref{deggrowth} to \re{equs}, and infer that $v(y)/M^{j-1} > y_n$ in $Q_1$, which implies that $u(x)/M^{j-1}> x_n $ in $Q$.

If $\rho< \rho^\prime/4$ (i.e. $Q$ is ``far" from the bottom boundary) we use the scaled form of the  interior weak Harnack estimate, Theorem \ref{kswh} (for a full statement see Corollary 4.8 in \cite{KS2}): for each $t<1/2$  and a solution $w$ of $\Lmin[w]\le  g(y)$, $w\ge0$ in $Q_{2t}(e)$, $g\in L^r(Q_{2t}(e))$, $r>p_0$,
\begin{equation}\label{whii}
\inf_{Q_t(e)} w\ge c{\left(t^{-n}\int_{Q_t(e)}w^\varepsilon\right)}^{{1}/{\varepsilon}}- Ct^{2-n/r}\|g\|_{L^r(Q_{2t}(e))}.
\end{equation}
Now $Q_{\rho/\rho^\prime}(e)\subset\{y_n>1/4\}$, so \eqref{tukk} implies that
\begin{equation}\label{whii2}
|\{v/M^j > 1/4\}\cap Q_{\rho/\rho^\prime}(e)| \ge (1-\mu)  |Q_{\rho/\rho^\prime}(e)|
\end{equation}
and we can apply the scaled WHI \re{whii} (with $r=n$, $t=\rho/\rho^\prime<1/4$ and $g=\rho^\prime \tilde f$, $w=v$) to the inequality \re{equs}, noting also that  $Q_{2\rho}(x_0) \subset \{x_n<d_0\}$ (since $x_{0,n}+\rho<3\rho^\prime/4<15d_0/16$), and
\begin{equation}\label{whii3}
\|\rho^\prime \tilde f\|_{L^n(Q_{2\rho/\rho^\prime}(e))} = \|f\|_{L^n(Q_{2\rho}(x_0))}\le C(n,q)\|f\|_{L^\prime(Q_2)}\le C\delta_0.
\end{equation}
Thus by \re{whii} and \re{whii2} we  conclude, by increasing $M$ (i.e. diminishing $k$) and diminishing $\delta_0$ if necessary, that in $Q_{\rho/\rho^\prime}(e)$
\begin{equation}\label{whii4}
v\ge c_1 M^j - C_1 \delta_0 \ge M^{j-1/2} - 1 \ge  M^{j-1}> M^{j-1}y_n.
\end{equation}
This means that $u/x_n > M^{j-1}$ in $Q$. \smallskip

{\it Case 2.2.} $\rho< d_0/4$ and $Q\cap\{x_n<d_0/2\}=\emptyset$.
 We use the same scaling, but now the zoom constant $\rho^\prime$ stays away from zero, in fact $\rho^\prime>d_0$. We apply in exactly the same manner \re{whii} with $r=\min\{p,q\}$, again with $t=\rho/\rho^\prime<1/2$ and $w=v$, to \re{equs} in $Q_{2\rho/\rho^\prime}(e)$. So if  $ n\le p $ the inequalities \re{whii3} and \re{whii4} hold (with $C$ in  \re{whii3} depending also on $p$), while if $p<n$ we get
 that
 $$
 \|\rho^\prime \tilde f\|_{L^p(Q_{2\rho/\rho^\prime}(e))} = (\rho^\prime)^{1-n/p} \|f\|_{L^p(Q_{2\rho}(x_0))}\le (\rho^\prime)^{1-n/p}\delta_0,
 $$
so by  \re{whii}
$$
v\ge c_1 M^j - C_1 (\rho/\rho^\prime)^{2-n/p} (\rho^\prime)^{1-n/p} \delta_0 \ge  c_1 M^j - C_1 (1/4)^{2-n/p} d_0^{1-n/p} \delta_0,
$$
in $Q_{\rho/\rho^\prime}(e)$, and we finish as in the previous case, by increasing $M$ and diminishing $\delta_0$, if necessary. \hfill $\Box$
  \smallskip

It remains to prove the growth lemma, Lemma \ref{deggrowth}. It is actually more convenient for notations to state it in the following (equivalent) form.

\begin{thm}[growth lemma]\label{deggrowth1} Let $\Omega$ be a smooth bounded domain. For any $\nu>0$, there exist $k,a>0$ depending on $n$, $p$, $q$, $\lambda$, $\Lambda$, $\nu$, $d_0$, $\|b\|_{L^q(\Omega)}$, $\Omega$, such that if $u$ is a viscosity solution of
$$\Lmin[u]\le  f(x), \quad f,u\ge0\mbox{ in }\Omega, \qquad\mbox{and}\qquad \|f\|_{\lpr(\Omega)} \le a,$$ and we have
$$
|\{x\in \Omega \::\: u(x) \ge d(x)\}|\ge \nu |\Omega|.
$$
then $u> kd$ in $\Omega$.
\end{thm}

We recall that $d(x)=\mathrm{dist}(x,\partial\Omega)$ and that $d_0$ is the width of a neighbourhood of $\partial \Omega$ where $f\in L^q$.
We get Lemma \ref{deggrowth} by applying Theorem \ref{deggrowth1} for some fixed smooth domain $\Omega$ such that $Q_{3/2}\subset \Omega\subset Q_2$.\smallskip

\noindent {\it Proof of Theorem \ref{deggrowth1}}. Take $d_1\in (0,d_0)$ for which the set $\Omega_{\delta}:=\{x\in\Omega\,:\, \mathrm{dist}(x,\partial\Omega)<\delta\}$ be smooth and has measure such that $|\Omega_{\delta}|\le\nu/2|\Omega|$, for all $0<\delta\le d_1$. Then obviously, setting $S_{\delta}=\Omega\setminus\Omega_{\delta}$,
$$
|\{u\ge \delta\}\cap S_\delta|\ge|\{u\ge d\}\cap S_\delta|\ge \frac{\nu}{2}|S_\delta|.
$$
We can apply the interior weak Harnack inequality (or the interior growth lemma, see \cite{Saf80} and \cite{SirARMA} for a direct proof of that result)
 and deduce that for each $\delta\in(0,d_1)$ we can find $k_\delta^\prime>0$ such that $$u\ge k_\delta^\prime\qquad \mbox{in }\;S_\delta.$$

We introduce the auxiliary function
$$
v_\delta = \frac{1}{\delta}\,d^2 + d.
$$
It is not difficult to compute that $\mm(D^2v_\delta)\ge 0$ in $\Omega_\delta$ if $\delta\in (0,d_2)$, provided  $d_2<d_1$ is chosen small enough, depending only on $n$,  $\lambda$, $\Lambda$ and the curvature of $\partial\Omega$ (a full computation can be found for instance in \cite{QS}, page~130). So we have
$$\Lmin[v_\delta]\ge - b|Dv_\delta|\ge -b\left(\frac{2d}{\delta} |Dd|+|Dd|\right)  \ge - 3b \qquad\mbox{in }\; \Omega_\delta,
$$
Set $$k_\delta=\frac{2\delta}{k_\delta^\prime}.$$ We have $v_\delta=2\delta\le k_\delta u$ on $\partial S_\delta$, $v_\delta=0\le k_\delta u$ on $\partial \Omega$, and
\begin{equation}\label{lastus}
\Lpl[v_\delta-k_\delta u]\ge \Lmin[v_\delta] - k_\delta\Lmin[u] \ge -3b - k_\delta f\qquad\mbox{in }\; \Omega_\delta.
\end{equation}
 We apply the ABP inequality, Theorem \ref{ksabp}  with $p=q$ and $R_0=2\delta$, to the inequality \re{lastus}. We deduce that
$$
v_\delta-k_\delta u\le {C} \delta^{1+\alpha}( \|b\|_{L^q(\Omega_\delta)} + k_\delta\|f\|_{L^q(\Omega_\delta)})\qquad\mbox{in }\;\Omega_{\delta},
$$
where $\alpha=1-n/q>0$  and ${C}$ depends  on $n,q,\lambda, \Lambda$, $d_0$, $\|b\|_{L^q}$, $\Omega$. Therefore, by the Lipschitz bound at $\partial \Omega$ applied to \re{lastus} (see Theorem \ref{lip}.2.)
\begin{eqnarray*}
v_\delta-k_\delta u&\le& \frac{C}{\delta}\left( \sup_{\Omega_\delta}(v_\delta-k_\delta u) + \delta^{2-n/q}\|3b +k_\delta f\|_{L^q(\Omega_\delta)}\right) \,d \\
&\le& {C} \delta^{\alpha}( \|b\|_{L^q(\Omega_\delta)} + k_\delta\|f\|_{L^q(\Omega_\delta)})\, d\qquad\mbox{in }\;\Omega_{\delta/2}.
\end{eqnarray*}
Note that we used the rescaled version of Theorem \ref{lip}.2. in a domain with width~$\delta$. Thus, since $v_\delta\ge d$, we see that
$$k_\delta u\ge (1 - {C} \delta^{\alpha}  \|b\|_{L^q(\Omega_\delta)}  - C \delta^{\alpha}k_\delta\|f\|_{L^q(\Omega_\delta)}) )\,d\qquad\mbox{in }\;\Omega_{\delta/2}.
$$
We now fix $\delta_0\in (0,d_2)$ small enough so that ${C} \delta_0^{\alpha}  \|b\|_{L^q(\Omega)}<1/4$ and then fix $a>0$ small enough so that $C \delta_0^{\alpha}k_{\delta_0} a<1/4$.

Thus $u\ge 1/(2k_{\delta_0})\, d$ in $\Omega_{\delta_0/2}$. On the other hand $u\ge (k_{\delta_0/2}^\prime/\mathrm{diam}(\Omega))\, d$ in $S_{\delta_0/2}$.
This ends the proof, setting
$$
k= \min\left\{\frac{1}{2k_{\delta_0}}, \frac{k_{\delta_0/2}^\prime}{\mathrm{diam}(\Omega)}\right\}=
\frac{k_{\delta_0/2}^\prime}{\mathrm{diam}(\Omega)},
$$
if $\delta_0$ is set to be smaller than $\mathrm{diam}(\Omega)/4$.\hfill $\Box$

\ed

{Some history}

\textsc{Kazdan---Kramer} 1975,
various partial results\bigskip

\textsc{Boccardo---Murat---Puel} 1980-1985\\
succeeded to obtain a full solvability result for the general class of divergence form operators, in the case when
$$
c(x)\le -c_0<0
$$\bigskip

\textsc{Ferone---Murat} 2000\\
solvability for the case
$$c\equiv0$$
and observation that only possible if $|Mh|$ small. \bigskip

\textsc{Barles-Murat, B.-Porretta}... 1995-2000 --- uniqueness. \bigskip

{\small
Related results by dall'Aglio, Giachetti and Puel; Maderna, Pagani and Salsa; Grenon, Murat and Porretta; Abdellaoui, dall'Aglio and Peral;  Abdel Hamid and Bidaut-V\'eron; Boccardo, Gallouet, Murat...}
\end{frame}

\documentclass[a4paper]{article}

\newtheorem{thm}{Theorem}
\newenvironment{thmbis}[1]
  {\renewcommand{\thethm}{\ref{#1}$'$}%
   \addtocounter{thm}{-1}%
   \begin{thm}}
  {\end{thm}}

\begin{document}
\begin{thm}
$1+1=2$
\end{thm}
\begin{thm}\label{comm}
$a+b=b+a$
\end{thm}

\begin{thm}
$0\ne0$
\end{thm}
\begin{thmbis}{comm}
$x+y=y+x$
\end{thmbis}

\begin{thebibliography}{99}\small

\bibitem{ARV} M.E. Amendola, L. Rossi, A. Vitolo,
Harnack inequalities and ABP estimates for nonlinear second-order elliptic equations in unbounded domains.
Abstr. Appl. Anal. 2008 (2008), Art. ID 178534.

\bibitem{AS} S.N. Armstrong, B. Sirakov,
 Nonexistence of positive supersolutions of elliptic equations via the maximum principle. Comm. Part. Diff. Eq. 36 (2011), 2011-2047.

\bibitem{ASS}
S. Armstrong, B. Sirakov, C. Smart,  Singular solutions of fully nonlinear elliptic equations and applications. Arch. Rat. Mech. Anal. 205 (2) (2012), 345-394.

\bibitem{B}   P.  Bauman,  Positive solutions of elliptic equations in nondivergence form and their adjoints. Ark. Mat. 22 (2) (1984),  153-173.

\bibitem{BNV} H. Berestycki, L. Nirenberg, S.R.S. Varadhan, {
The principal eigenvalue and maximum principle for second order
elliptic operators in general domains},  Comm. Pure Appl. Math. {
47}(1) (1994),  47-92.

\bibitem{BEM} J. Braga,  M. Ederson, D. Moreira,  Inhomogeneous Hopf-Oleinik Lemma and Applications. Part I: Regularity of the Normal Mapping, in preparation.

\bibitem{BMW} J. Braga,  D. Moreira, L. Wang, Inhomogeneous Hopf-Oleinik lemma and Krylov's boundary gradient estimates, arXiv:1608.02352.

\bibitem{BC} H. Brezis, X. Cabr\'e,
Some simple nonlinear PDE's without solutions.
Boll. Unione Mat. Ital.  8 (1) (1998), 223-262.

\bibitem{Cab} X. Cabr\'e, On the Alexandroff-Bakelʹman-Pucci estimate and the reversed H\"older inequality for solutions of elliptic and parabolic equations. Comm. Pure Appl. Math. 48 (5) (1995), 539-570.

\bibitem{CC}
L. Caffarelli, X. Cabr\'e, Fully nonlinear elliptic equations. American Mathematical Society Colloquium Publications, 43, Providence, (1995).

\bibitem{CCKS}
L. Caffarelli, M. Crandall, M. Kocan, A. Swiech, On viscosity solutions of fully nonlinear equations with measurable ingredients,
Comm. Pure  Appl. Math.
{\bf 49}(4) (1996),  365-398.

\bibitem{CCCC}   L. Caffarelli, E. Fabes, S. Mortola, S. Salsa,  Boundary behavior of nonnegative solutions of elliptic operators in divergence form. Indiana Univ. Math. J. 30 (4) (1981),  621-640.

\bibitem{CLN} L. Caffarelli, Y. Li, L. Nirenberg, Some remarks on singular solutions of nonlinear elliptic equations III: viscosity solutions including parabolic operators. Comm. Pure Appl. Math. 66 (1) (2013),  109-143.

    \bibitem{CIL}
 M. Crandall, H. Ishii and P.L. Lions. User's guide to viscosity
solutions of second order partial differential equations.
Bulletin of the AMS, {\bf 27}(1) (1992),  1-67.

 \bibitem{GS}   M. G. Crandall, A. Swiech, A note on generalized maximum principles for elliptic and parabolic PDE,  in Evolution Equations, Proceedings in honor of J. A. Goldstein's 60th birthday (Goldstein, Nagel and Romanelli eds.), Lecture notes in pure and applied mathematics, vol. 234 (2003), 121-127.

 \bibitem{E}    L. Escauriaza, $W^{2,n}$ a priori estimates for solutions to fully non-linear equations, Indiana Univ.
Math. J. 42 (1993), 413-423.

 \bibitem{FS} E. B. Fabes, D. W. Stroock, The $L^p$-integrability of Green’s functions and fundamental
solutions for elliptic and parabolic equations, Duke Math. J. 51 (1984), 997-1016.



 \bibitem{GT} D. Gilbarg, N.S. Trudinger, { Elliptic Partial
Differential Equations of Second Order}, 2nd edition, Springer
Verlag.

\bibitem{KS1}    S. Koike,  A. Swiech,  Existence of strong solutions of Pucci extremal equations with superlinear growth in Du. J. Fixed Point Theory Appl. 5 (2009), no. 2, 291-304.

\bibitem{KS2}    S.  Koike, A. Swiech, Weak Harnack inequality for fully nonlinear uniformly elliptic PDE with unbounded ingredients. J. Math. Soc. Japan 61 (3) (2009),  723-755.

\bibitem{KS3}    S.  Koike, A. Swiech, Local maximum principle for Lp-viscosity solutions of fully nonlinear elliptic PDEs with unbounded coefficients. Comm. Pure Appl. Anal. 11 (5) (2012),  1897-1910.

\bibitem{K0}    N.V.  Krylov, Boundedly inhomogeneous elliptic and parabolic equations in a domain.
Izv. Akad. Nauk SSSR Ser. Mat. 47(1) (1983), 75-108.

\bibitem{Kbook} N.V. Krylov, Nonlinear elliptic and parabolic
equations of second order. Mathematics and its Applications,
Reidel, 1987.

\bibitem{K} N.V. Krylov,  Some $L^p$-estimates for elliptic and parabolic operators with measurable coefficients. Discr. Cont. Dyn. Syst. B 17 (6) (2012), 2073-2090.

       \bibitem{QS} A. Quaas, B. Sirakov, Principal eigenvalues and the Dirichlet poblem for fully nonlinear elliptic operators. Adv. Math. 218 (2008), 105-135.

\bibitem{Saf80} M. V. Safonov,  Harnack's inequality for elliptic equations and H\"older property of their solutions. (Russian)  Zap. Nauchn. Sem. Leningrad. Otdel. Mat. Inst. Steklov. (LOMI) 96 (1980), 272-287. English transl. in J. Soviet Math, 21 (5) (1983), 851-863.

\bibitem{Saf}   M. V.   Safonov,  Non-divergence elliptic equations of second order with unbounded drift. Nonlinear partial differential equations and related topics, 211-232, Amer. Math. Soc. Transl. Ser. 2, 229, Amer. Math. Soc., Providence, RI, 2010.

\bibitem{SirARMA} B. Sirakov, Solvability of uniformly elliptic fully nonlinear PDE, Arch. Ration. Mech. Anal. 195 (2) (2010), 579-607.

 \bibitem{S2} B. Sirakov, Uniform bounds via regularity estimates for elliptic PDE with critical growth in the gradient, arXiv:1509.04495

     \bibitem{W} N. Winter, $W^{2,p}$ and $W^{1,p}$-estimates at
the boundary for solutions of fully nonlinear,
uniformly elliptic equations,
Z. Anal. Adwend. (Journal for Analysis and its Applications) 28(2) (2009), 129-164.
\end{thebibliography}
\end{document}